\documentclass[a4paper, 12pt]{article}
\usepackage[koi8-r,cp1251]{inputenc}
\usepackage[english,russian]{babel}
\usepackage{amsfonts,amssymb,amsmath, hyperref}
\usepackage[final]{epsfig}
\usepackage{graphicx}

\begin{document}

\begin{center}
THE PALEY-WIENER-GELFAND TAUBERIAN THEOREM FOR  SEMIGROUPS WITH INVARIANT MEASURE\\
 Mirotin A. R.\\
 amirotin@yandex.ru\\
 \end{center}
 
 \
 
\small{The theorem is proved that generalizes the Gelfand generalization of the Paley-Wiener tauberian
 theorem to general abelian topological semigroups with invariant measure. More precisely,  it is shown that if $w$ is a weight on an Abelian invariant measure  semigroup $S,$  $fw\in
L^\infty(S),\ K\in L^1(w)$, and  $\widetilde
K(\psi)\ne -1$ for all  $\psi\in \widetilde S_w,$ then $(f+K\ast
f)w\to 0$ on $\infty$ implies $fw\to 0$ on $\infty.$
Several corollaries
 of this theorem are given.}

\

\begin{center}
{\bf    ТАУБЕРОВА ТЕОРЕМА ПЭЛИ-ВИНЕРА-ГЕЛЬФАНДА ДЛЯ ПОЛУГРУПП С ИНВАРИАНТНОЙ МЕРОЙ\\}
\end{center}

\

\small{Доказана теорема, переносящая обобщение Гельфанда тауберовой теоремы Пэли-Винера на
общие абелевы топологические полугруппы с инвариантной мерой. Приведено несколько следствий
этой теоремы.}

\

\begin{center}
{\bf  1. Ведение и предварительные сведения }
\end{center}

Как известно, теоремы, в которых асимптотическое поведение
некоторого усреднения функции выводится из асимптотического
поведения самой этой функции, называются {\it абелевскими}, а
теоремы, в которых утверждается обратное, носят название {\it
тауберовых}.  Как правило, справедливость тауберовой теоремы
требует наличия некоторого дополнительного (по сравнению с
соответствующей абелевской теоремой) условия, которое принято
называть {\it тауберовым.} Терминология восходит к классическим
результатам Абеля и Таубера, в которых речь шла о рядах. Первоначальным стимулом к развитию этого направления в анализе послужили задачи теории чисел, в дальнейшем появились и другие приложения, в частности, к теории операторов
(см. \cite{Kor}).

Р. Пэли и Н. Винеру \cite{PW}  принадлежит  теорема тауберовского типа,
утверждающая, что {\it для функции $k$ из $L^1(\Bbb{R}_+)$, такой,
что её преобразование Лапласа $\tilde{k}(z)$ не принимает значения
-1 при ${\rm Re}z\geq 0$, и функции  $f$ из $L^\infty(\Bbb{R}_+)$
 соотношение $f(x)+f\ast k(x)\to s$  влечет соотношение}
$f(x)\to s/(1+\tilde{k}(1))\quad (x\to +\infty).$ Ими же показано,
что тауберово условие здесь существенно. В качестве частного
случая получается следующая теорема Мерсера: {\it если
$0<\alpha<1$ и если $\alpha x_n+(1-\alpha)n^{-1}\sum_{k=1}^n
x_k\to x$, то} $x_n\to x$. Теорема Пэли-Винера на
абелевы полугруппы с инвариантной мерой  была перенесена в \cite{Trudy}.

 И. М. Гельфанд  \cite{GRS}  обобщил
результат Пэли и Винера на случай весовой алгебры $L^1(\Bbb{R}_+, w)$. Целью настоящей работы является
перенесение тауберовой теоремы Пэли-Винера-Гельфанда на общие
абелевы полугруппы с инвариантной мерой и получение некоторых её
многомерных и дискретных аналогов.

Далее $S$  будет обозначать абелеву топологическую полугруппу,
${\cal B}(S)$ -- $\sigma$-алгебру  борелевских подмножеств $S$.
Полугруппа $S$ будет предполагаться наделенной инвариантной мерой,
а именно, ненулевой
положительной  мерой Радона $\mu$, такой, что $\mu(sB)=\mu(B)$ для
любых $B\in{\cal B}(S)$, обладающих тем свойством, что $sB\in{\cal
B}(S)\quad (s\in S)$. Мера $\mu$ будет  предполагаться также
локально ограниченной в том смысле, что каждая точка обладает
окрестностью конечной меры. Простейшие примеры таких полугрупп мы получим,
рассматривая борелевские подполугруппы локально компактных абелевых групп с
суженной на них мерой Хаара. Другим примером является плоскость Немыцкого с
обычным сложением и плоской мерой Лебега. Так как $S$ есть полугруппа с
инвариантной мерой в смысле \cite{SF}, мы будем
 использовать построенную там структурную теорию таких полугрупп
 (см. также \cite[глава 1]{garman}).
В частности, в \cite{SF} было показано, что существует локально компактная абелева группа $G$,
ненулевой гомоморфизм $\tau$ полугруппы $S$ на подполугруппу $\Sigma$ группы $G$ и идеал $S_0$
полугруппы $S$, такие, что 1) сужение $\tau_0$ гомоморфизма $\tau$ на $S_0$ есть топологическое вложение; 2) идеал $\tau(S_0)$ полугруппы $\Sigma$ открыт в $G$ и порождает $G$; 3) образ $\tau_0(\mu)$ и сужение меры Хаара  $\nu$ группы $G$ на  $\tau(S_0)$  совпадают; 4) $\mu(C)=\nu(\tau(C))$ для любых компактных $C\subset S$. Отсюда следует, что
для любой функции $\varphi\in L^1(G)$ имеет место равенство (\cite[следствие 1.2.5]{garman})
\begin{equation}\label{125}
\int\limits_\Sigma\varphi d\nu=\int\limits_S\varphi\circ\tau d\mu.
\end{equation}

Будем предполагать дополнительно, что $\Sigma=\tau(S)$ есть
борелевское подмножество группы $G$ (хотя часть результатов
справедлива и без этого предположения). Все функции на $\Sigma$ считаем доопределенными нулем на $G\setminus \Sigma$.

Нам также понадобятся некоторые свойства алгебр $L^1(w)$ и
$L^\infty(S)$, установленные в \cite{garman}.

{\bf Определение 1.} \textit{Полухарактером} полугруппы $S$ будем называть непрерывный гомоморфизм этой полугруппы в мультипликативную полугруппу $\mathbb{C}$, отличный от тождественно нулевого. Через $\widetilde S$  ($\widetilde S_+$)  мы обозначаем пространство
всех не локально $\mu$-пренебрежимых (соответственно всех не локально $\mu$-пренебрежимых неотрицательных) полухарактеров полугруппы
$S$, наделенное компактно-открытой топологией.

 {\bf Определение 2.} Неотрицательная измеримая функция $w$
на $S$ называется {\it весом},  если она локально
ограничена (т. е. ограничена в достаточно малой окрестности любой
точки $s\in S$), и справедливо неравенство
$$
w(st)\leq w(s)w(t)\quad(s,t\in S).
$$

Далее будем считать, что вес $w$ факторизуем, т. е. $w=w'\circ \tau$, где  $w'$ --- вес на $\Sigma$. (В силу следствия 1.2.2 из \cite{garman} каждый вес факторизуем, если $S$  --- полугруппа с сокращениями.)

Всюду ниже через $\widetilde S_w$  обозначается множество
всех $\gamma\in\widetilde S$, для которых $|\gamma|\leq w$. Если $w=\rho\in \widetilde  S_+$ --- неотрицательный полухарактер, то каждый полухарактер из $\widetilde S_\rho$ имеет вид $\psi=\rho\xi$, где  $\xi\in \widetilde  S$ --- ограниченный полухарактер.

Пространство $L^1(S,w\mu)$ будем  для краткости обозначать
$L^1(w)$.  Оно является коммутативной банаховой алгеброй со сверткой в качестве умножения, стандартная норма в ней будет обозначаться
$\|.\|_{w}$.

{\bf Определение 3}  \cite{VUZ88}. Пусть $k\in L^1(w)$. Функцию $\tilde k$
на $\widetilde S_w$, определенную равенством
$$
\tilde k(\gamma)=\int\limits_S k\overline\gamma d\mu,
$$
будем называть {\it преобразованием Лапласа}
 функции $k$.

Следующая теорема, доказанная в \cite{garman},  (см. там теорему 3.2.8)  является  обобщением результатов Гельфанда и Эдвардса (\cite[с. 370 -- 372]{Edv}). В частном случае $w\in \widetilde S_+$ она была опубликована в \cite{VUZ88}.

{\bf Теорема 1.} {\it Комплексные гомоморфизмы  алгебры
$L^1(w)$ есть в точности отображения вида  $f\mapsto\tilde
f(\gamma)$, где $\gamma\in\widetilde S_w$. Гельфандов спектр
алгебры $L^1(w)$ можно отождествить с $\widetilde S_w$, и при этом отождествлении
преобразование Гельфанда для $ L^1(w)$ есть преобразование
Лапласа.}

\begin{center}
{\bf  2. Одна теорема абелевского типа }
\end{center}

Всюду ниже вес $w$ считается положительным.

 {\bf  Определение 4.} Измеримую функцию $f$ на $S$ вида $f=f'\circ \tau$,  где $f'$ --- измеримая функция  на $\Sigma$, назовем \textit{факторизуемой}.

 {\bf  Определение 5.}  Для факторизуемой функции $f$ на $S$,  $f=f'\circ \tau$, определим \textit{сдвиг} на элемент $s\in S$ равенством
$$
f_s(t)=f'(\tau(s)^{-1}\tau(t)),\ t \in S.
$$

{\bf  Определение 6.} Пусть функции $f$ и $k$ определены на $S$, причем $f$ факторизуема.
Определим  \textit{свертку} функций $f$ и $k$ равенством
$$
  (f\ast k)(t)= \int\limits_S f_s(t)k(s)d\mu(s),\ t \in S,
$$
если интеграл в правой части существует для $\mu$-почти всех $t$. В этом случае
 функции $f$ и $k$ будем называть \textit{свертываемыми}.

Для получения основных результатов нам понадобится ряд лемм. Часть из них опубликована в \cite{garman}, мы приводим доказательства ввиду труднодоступности этого источника.

{\bf Лемма 1}. {\it Отображение $\tau':f'\mapsto f'\circ\tau$ есть изоморфизм
сверточных алгебр $L^1(w')$ и  $L^1(w)$. В частности, $L^1(w)$
есть коммутативная банахова алгебра.}

Доказательство. Достаточно убедиться в справедливости равенства
$k\ast g=(k'\ast g')\circ\tau$, а это сразу вытекает из формулы (\ref{125}).

{\bf Лемма 2}. {\it Отображение $h'\mapsto h'\circ\tau$ есть
изометрический изоморфизм банаховой алгебры $L^\infty(\Sigma)$ на
банахову алгебру}    $L^\infty(S)$.

Доказательство. Из леммы 1 следует, что отображение $f\mapsto
f\circ\tau'$ есть изометрический изоморфизм  сопряженных
пространств $L^1(S)'$ и $L^1(\Sigma)'$. В силу общего вида
линейных функционалов в пространствах $L^1$ отсюда выводим, что
отображение, которое каждой функции $h'\in L^\infty(\Sigma)$
ставит в соответствие функцию $h\in L^\infty(S)$, удовлетворяющую
при всех $\varphi\in L^1(\Sigma)$  условию
$$
\int\limits_\Sigma h'\varphi d\nu=\int\limits_S h(\tau'\varphi)d\mu,
$$
\noindent есть  изометрический изоморфизм пространства
$L^\infty(\Sigma)$ на пространство  $L^\infty(S)$ (это
отображение, обратное к $f\mapsto f\circ\tau'$). Но по формуле
(\ref{125}) при всех $\varphi\in L^1(\Sigma)$
$$
\int\limits_\Sigma h'\varphi d\nu= \int\limits_S
(h'\circ\tau)(\varphi\circ\tau) d\mu.
$$
\noindent Таким образом, с учетом леммы 1 при всех  $g\in
L^1(S)$ справедливо равенство
$$
\int\limits_S (h'\circ\tau)gd\mu= \int\limits_S hgd\mu,
$$
\noindent а потому $h=h'\circ\tau$.

{\bf Лемма  3.}  {\it Пусть функция $f$ на полугруппе $S$
измерима, $w$ ---   вес на $S$, $fw\in L^\infty(S)$, $k\in L^1(w)$.  Тогда функции $f$ и $k$ свертываемы, и}
$$|f\ast k|w\leq (|f|w)\ast (|k|w).$$

Доказательство. Функция $fw$ факторизуема в силу леммы 2. Кроме того, из легко проверяемого неравенства $\|(fw)_t\|_\infty\leq\|fw\|_\infty$ следует, что $fw$ и $kw$  свертываемы. В силу факторизуемости и положительности веса, факторизуема и сама функция $f$.  Используя введенные выше обозначения, имеем
$$
|f_s(t)|w(t)=
|f'(\tau(t)\tau(s)^{-1})|w'(\tau(t)\tau(s)^{-1}\tau(s))\leq
$$
$$
\leq |f'(\tau(t)\tau(s)^{-1})|w'(\tau(t)\tau(s)^{-1})w'(\tau(s))
$$
(обе части последнего неравенства равны нулю при $\tau(t)
\tau(s)^{-1}\notin\Sigma$). Таким образом, $|f_s(t)|w(t)\leq
(|f|w)_s(t)w(s)$. Следовательно,
$$
\int\limits_S|f_s(t)||k(s)|d\mu(s)w(t)\leq
\int\limits_S(|f|w)_s(t)||k(s)|w(s)d\mu(s)=(|f|w)\ast (|k|w)(t),
$$
что доказывает лемму.

{\bf  Лемма 4}. а) {\it   Для любых  $f\in L^\infty(S),
\varphi\in L^1(S)$ и  $t\in S$ справедливо равенство}
$$
 \int\limits_S f_t(s)\varphi(s) d\mu(s)=\int\limits_S f(s) \varphi(ts)d\mu(s).
$$
б) {\it Для любых  $k\in L^1(w), h\in L^\infty(S)$ и  $t\in S$
справедливо равенство}
$$
 \int\limits_S k_t(s)h(s)w(s) d\mu(s)=\int\limits_S k(s) h(ts)w(ts)d\mu(s).
$$

Доказательство. а). Пусть  $f'\in  L^\infty(\Sigma), \varphi'\in
L^1(\Sigma)$ таковы, что $f=f'\circ\tau, \varphi=\varphi'\circ\tau$.
Дважды используя формулу (\ref{125}), имеем
$$
\int\limits_S\!f_t\varphi d\mu\!=\!\int\limits_S\!
f'_{\tau(t)}\!(\tau(s))\!\varphi'\!(\tau(s)\!)d\mu(s)\!=\!
\int\limits_\Sigma\!f'_{\tau(t)}(x)\!\varphi'(x)d\nu(x)\!=
$$
$$
\hspace{7mm}=\!\int\limits_\Sigma f'(z)\varphi'(\tau(t)z)d\nu(z) =
\int\limits_S f'(\tau(s)) \varphi'(\tau(ts))d\mu(s)=
$$
$$
= \int\limits_S f(s) \varphi(ts)d\mu(s).\hspace{35mm}
$$

Доказательство пункта б) проводится аналогично.

Ниже через $1_E$ обозначается индикатор множества $E\subset S$.

{\bf  Лемма 5}. {\it {\rm а)} Для любых $E\in {\cal B}(S),\
t\in S$
 таких, что $tE\in {\cal B}(S)$, справедливо равенство
 $(1_E)_t =1_{tE}$   {\rm (}в смысле алгебры}   $L^\infty(S))$;

б) $\|f_t\|_\infty\leq \|f\|_\infty \mbox{ для любых}  \  f\in
L^\infty(S), t\in S.$

Доказательство. а). Для любого компактного $D\subset S$ имеем в
силу леммы 4
$$
   \int\limits_D(1_E)_td\mu=\int\limits_S 1_D(1_E)_td\mu= \int\limits_S
1_D(ts)1_E(s)d\mu(s)=\hspace{2cm}
$$
$$
\hspace{6mm}=\int\limits_S 1_{t^{-1}D}(s)1_E(s)d\mu(s)=
\mu(t^{-1}D\cap E),
$$
\noindent где $t^{-1}D:=\{s\in S:ts\in D\}$ -- борелевское
множество. Но, поскольку $t(t^{-1}D\cap E)=D\cap tE$, то
$$
\mu(t^{-1}D\cap E)=\mu(D\cap tE)=\int_D 1_{tE}d\mu.
$$
\noindent Следовательно, $(1_E)_t=1_{tE}\quad \mu$-п.в.

б). Это является непосредственным следствием определений.

{\bf  Лемма 6.} {\it Для  $f\in L^\infty(S),\ k\in L^1(S)$
справедливо равенство
$$(f\ast k)'=f'\ast k',
$$
\noindent где $f=f'\circ\tau, k=k'\circ\tau$ и }$f'\in
L^\infty(\Sigma),\ k'\in L^1(\Sigma)$.

Доказательство. Применяя формулу (\ref{125}), получаем
$$
(f\ast k)(t)= \int\limits_S f_s(t)k(s)d\mu(s)=\int\limits_S
f'_{\tau(s)}(\tau(t))k'(\tau(s))d\mu(s)=
$$
$$
\hspace{3mm}=\int\limits_\Sigma f'_x(\tau(t))k'(x)d\nu(x)=(f'\ast
k')(\tau(t)).
$$

{\bf  Лемма 7.}  {\it Пусть $f\in L^\infty(S),\ k, g\in
L^1(S)$. Тогда}

а) $f\ast k\in L^\infty(S)\cap C(S)$ , $\quad \|f\ast
k\|_\infty\leq\|f\|_\infty\|k\|_1$;

б) $(f\ast k)_t=f\ast k_t$;

в) $(f\ast k)\ast g=f\ast (k\ast g)$;

г)  $(f\ast k)(t)=\int_S f(s)k_s(t)d\mu(s)$.

Доказательство. а). В силу леммы 5 имеем
$$
  |(f\ast k)(t)|\leq  \int\limits_S |f_s(t)||k(s)|d\mu(s)\leq \|f\|_\infty
  \|k\|_1.
$$
Осталось установить непрерывность свертки (т. е. показать, что
соответствующий класс эквивалентности содержит непрерывную
функцию). Но известно (см., например, \cite[гл. 5, \S 20]{HiR1}),
что свертка $f'\ast k'$  $\nu$-почти всюду совпадает с некоторой
функцией $h'\in L^\infty(\Sigma)\cap C(\Sigma)$. Следовательно,
$f\ast k=h'\circ\tau$  $\mu$-п.в. (В противном случае $(f'\ast
k')\circ\tau\ne h'\circ\tau$ на некотором компакте $C\subset S$
положительной $\mu$-меры, а тогда    $f'\ast k'\ne h'$ на
$\tau(C)$, причем  $\nu(\tau(C))=\mu(C)>0$, --- противоречие).

б).  С учетом леммы 6 получаем, что
 $$(f\!\ast\! k)_t(s)\!=\!(f\!\ast\!
 k)'_{\tau(t)}\!(\!\tau(s)\!)\!=\!(\!f'\!\ast\!k'\!)_{\tau(t)}\!(\!\tau(s)\!)\!
=\!(f'\!\ast\!k'_{\tau(t)})\!(\!\tau(s)\!),
$$
\noindent причем по определению $k'_{\tau(t)}=(k_t)'$.

в). Это легко выводится из соответствующего свойства свертки на
группе $G$ с помощью леммы 6.

 г). Снова используя свойства
свертки на группе $G$ и лемму 6, имеем в силу формулы (\ref{125})
$$
\hspace{-12mm}(f\ast k)(t)=(f'\ast
k')(\tau(t))=\int\limits_{\Sigma} f'(x)k'_x(\tau(t))d\nu(x)=
$$
$$
\hspace{8mm}=\int\limits_S
f'(\tau(s))k'_{\tau(s)}(\tau(t))d\mu(s)=\int\limits_S
f(s)k_s(t)d\mu(s).
$$

 Будем говорить, что функция $f$ на
полугруппе $S$ {\it стремится к числу $c$ на }$\infty$ (и писать $f\to c$ на $\infty$), если для любого
$\epsilon>0$ вне некоторого компакта
$K_\epsilon\subset S$ выполняется неравенство
$|f(x)-c|<\epsilon$.

Теперь мы можем установить  следующее утверждение абелевского типа.

{\bf  Теорема 2}.  {\it  Если функция $F$ измерима на $S$,
$Fw\in L^\infty(S)$, $k\in L^1(w)$ и $Fw\to 0$ на $\infty$, то
существует $F\ast k$, и $(F\ast k)w\to 0$ на} $\infty$.

 Доказательство. Существование $F\ast k$ сразу следует из леммы 3.
 Положим $f_1:=Fw,\ k_1:=kw$. В силу той же леммы  можно считать, что
 $f_1\ne 0,\ k_1\ne 0$. Для любого $\epsilon
>0$  выберем такие  компактные $C_1,\ C_2\subset S$, что
$$
\int\limits_{S\setminus
C_1}|k_1|d\mu<\frac{\epsilon}{2\|f_1\|_\infty}, \mbox{ и }
|f_1(s)|<\frac{\epsilon}{2\|k_1\|_1} \mbox{ при } s\in S\setminus
C_2.
$$
\noindent Очевидно, что
$$
|f_1|\ast |k_1|(t)=\int\limits_{C_1}
|(f_1)_s(t)||k_1(s)|d\mu(s)+\int\limits_{S\setminus  C_1}
|(f_1)_s(t)||k_1(s)|d\mu(s).
$$
Пусть $t\notin C_1C_2$. Тогда в силу лемм 5  и  7
 имеем
$$
\hspace{-11mm}\int\limits_{C_1}
|(f_1)_s(t)||k_1(s)|d\mu(s)=\int\limits_{S}
|(f_1)_s(t)||k_1(s)|1_{C_1}(s)d\mu(s)=
$$
$$
=|f_1|\ast|k_11_{C_1}|(t)=\int\limits_{S}
|f_1(s)|(|k_1|1_{C_1})_s(t)d\mu(s)=
$$
$$
\!=\!\int\limits_{S\setminus C_2}\!
+\!\int\limits_{C_2}\!|f_1(s)||(k_1)_s(t)|1_{sC_1}(t)d\mu(s)\!
\leq\!\int\limits_{S\setminus C_2}\!|f_1(s)||(k_1)_s(t)|d\mu(s)\!+
$$
$$
 +\int\limits_{C_2}|f_1(s)||(k_1)_s(t)|1_{sC_1}(t)d\mu(s) \leq\frac{\epsilon}{2\|k_1\|_1} \int\limits_{S}
|(k_1)_s(t)|d\mu(s),
$$
\noindent поскольку $1_{sC_1}(t)=0$ при $t\notin C_1C_2, \ s\in
C_2$.

Далее,  имеем для некоторого $k'\in
L^1(\Sigma)$
$$
\int\limits_{S}|(k_1)_s(t)|d\mu(s)=\int\limits_{S}|k'(\tau(s)^{-1}\tau(t))|d\mu(s).
$$
\noindent Рассмотрим функцию $h(x)=|k'(x^{-1}\tau(t))|$ из
$L^1(G)$.  Учитывая инвариантность меры Хаара $\nu$ относительно
сдвигов и отражения, получаем
$$
\hspace{-14mm}\int\limits_{S}|(k_1)_s(t)|d\mu(s)=\int\limits_{S}h(\tau(s))d\mu(s)=
\int\limits_{\Sigma}h(x)d\nu(x)\leq
$$
$$
\hspace{14mm}\leq\int\limits_{G}|k'(x^{-1}\tau(t))|d\nu(x) =
\int\limits_{G}|k'(x)|d\nu(x)=\|k_1\|_1.
$$
Таким образом, $\int_{C_1} |(f_1)_s(t)||k_1(s)|d\mu(s)\leq
\epsilon/2.$  Но
$$
\int\limits_{S\setminus  C_1}
|(f_1)_s(t)||k_1(s)|d\mu(s)\leq\|f\|_\infty\int\limits_{S\setminus
C_1} |k_1(s)|d\mu(s)<\frac{\epsilon}{2}.
$$
Для завершения доказательства осталось заметить, что  $|F\ast k|w\leq |f_1|\ast |k_1|$ по лемме
3.

\begin{center}
{\bf  3. Тауберовы теоремы}
\end{center}

В условиях  предыдущей теоремы  $(F+F\ast k)w\to 0$  на $\infty$.
Оказывается, справедливо частичное обращение этого утверждения.
Теоремы 1 и 2   позволяют получить этот результат как
непосредственное следствие  гельфандовской теории коммутативных банаховых алгебр (ср. \cite[c. 124]{GRS}).

{\bf  Теорема 3}. {\it Пусть  $f$ измерима на $S$, $fw\in
L^\infty(S),\ K\in L^1(w)$, и  $\widetilde
K(\psi)\ne -1$  при всех $\psi\in \widetilde S_w$. Если $(f+K\ast
f)w\to 0$ на $\infty$, то и $fw\to 0$ на} $\infty$.

 Доказательство. Напомним, что по теореме 1
гельфандов спектр  алгебры $L^1(w)$ естественным образом
отождествляется с $\widetilde S_w$, и при этом преобразование
Гельфанда переходит в преобразование Лапласа. Обозначим через $V$
алгебру, получающуюся из $L^1(w)$ присоединением единицы $u$ (если
это необходимо), умножение в $V$ тоже обозначим $\ast$. Для
измеримой на $S$ функции $F$ и элемента $z=\lambda u+k\in V$
положим $z\ast F:=\lambda F+k\ast F$, при условии, что свертка
существует. Тогда, если свертки существуют, то $z_1\ast(z_2\ast
F)=(z_1\ast z_2)\ast F$ при $z_1, z_2\in V$. Тауберово условие
означает, что элемент $u+K$ обратим в $V$. Если мы положим
$F:=(u+K)\ast f$, то $Fw=fw+(K\ast f)w\in L^\infty(S)$ в силу лемм
3 и 7 а). Теперь в силу  теоремы 2
$((u+K)^{-1}\ast F)w\to 0$ на $\infty$, причем свертка существует,
и для завершения доказательства осталось заметить, что $(u+K)^{-1}\ast
F=((u+K)^{-1}\ast(u+K))\ast f=f$.

Если учесть, что  непрерывные полухарактеры полугруппы
$\Bbb{R}_+^n$ имеют вид $x\mapsto e^{-x\cdot z}\ (z\in\Bbb{C}^n$), то из теоремы 3 вытекает следующая тауберова теорема для
$n$-мерной усеченной свертки (ниже для $x\in \Bbb{R}_+^n$ положено
$[0,x]=(x-\Bbb{R}_+^n)\cap\Bbb{R}_+^n$).

 {\bf  Следствие 1.}  {\it Пусть $w$ --- положительный вес на аддитивной
 полугруппе
 $\Bbb{R}_+^n$,   $fw\in
L^\infty(\Bbb{R}_+^n),\ k\in L^1(w)$, причем
$$
\int\limits_{\Bbb{R}_+^n}k(x)e^{-x\cdot z}dx\ne -1
$$
при  $z\in\Bbb{C}^n$ таких, что $x\cdot{\rm Re}z\geq -\log w(x)$
 при всех $x\in \Bbb{R}_+^n.$
 \noindent  Если
 $$
 \left(f(x)+\int_{[0,x]}f(x-y)k(y)dy\right)w(x)\to 0,\ x\to +\infty,
 $$
 то
 $$
 f(x)w(x)\to 0,\ x\to +\infty.
 $$}

{\bf  Следствие 2.}  {\it Пусть $S$ -- дискретная полугруппа с
единицей $e$ и сокращениями, $w$ --- положительный вес на этой
 полугруппе, $fw\in L^\infty(S), \ k\in L^1(w)$,
причем $\tilde k(\psi)\ne 0$  при всех $\psi\in \widetilde S_w$.
Если $(f\ast k)w\to 0$  на $\infty$, то $fw\to 0$ на} $\infty$.

 Действительно, в этом случае алгебра $L^1(w) (=l_1(w))$ содержит
единицу $u=1_{\{e\}}$. Тогда $(f\ast k)w=(f+f\ast (k-u))w$, причем
$\widetilde{(k-u)}(\psi)\ne -1$ при всех $\psi\in \widetilde S_w.$

{\bf  Следствие 3.}  {\it Пусть  $w$ --- положительный вес на аддитивной полугруппе $\mathbb{Z}_+$,
 $fw\in l^\infty(\mathbb{Z}_+), \ k\in l^1(w)$,
причем $\sum_{n=0}^\infty k(n)z^n\ne 0$  при всех $z\in \mathbb{C}$, таких, что $|z|\leq \inf_n\sqrt[n]{w(n)}$.
Если $\sum_{i=0}^\infty f(i) k(n-i)w(n)\to 0$ при  $n\to\infty$, то $f(n)w(n)\to 0$ при  $n\to\infty$.}

Это сразу вытекает из следствия 2, если заметить, что полухарактеры полугруппы $\mathbb{Z}_+$
имеют вид $\psi(n)=z^n,\ z\in \mathbb{C}\ (0^0:=1)$.

С помощью следствия 2 мы сейчас получим тауберову теорему для
свертки Дирихле.

{\bf  Теорема 4.}  {\it  Пусть $\rho$ --- положительный полухарактер
 мультипликативной полугруппы $\Bbb{N}^*$ натуральных чисел,
$f\rho\in L^\infty(\Bbb{N}^*)$,  $k\in L^1(\rho)$, причем сумма
соответствующего ряда Дирихле ограничена снизу в правой
полуплоскости:
$$
 \left|\sum\limits_{n=1}^{\infty}\frac{k(n)\rho(n)}{n^s}\right|\geq C>0\quad  (s\in
\Bbb{C}, \quad {\rm Re}s\geq 0).
$$

Если $\sum_{d|n}f(d)k(n/d)\rho(n)\to 0$, то $f(n)\rho(n)\to 0$ при}
$n\to\infty.$

  Доказательство сразу вытекает  из следствия 2 и следующей
  леммы, обобщающей теорему 1 из \cite{HW}, где рассмотрен случай $\rho=1$.

{\bf  Лемма 8.}  {\it  Пусть $\rho$ --- положительный полухарактер
 полугруппы $\Bbb{N}^*$,  $k\in
L^1(\rho)$. Если сумма  ряда Дирихле ограничена снизу в
правой полуплоскости  как в теореме 4, то $\widetilde k(\psi)\ne 0$  при всех} $\psi\in
\widetilde{\Bbb{N}^*}_\rho$.

  Доказательство. Так как $S=\Bbb{N}^*$ -- свободная полугруппа, образующими
которой являются простые числа, то общий вид полухарактера в
$\Bbb{N}^*$ есть $\psi(n)=z^{\alpha(n)}:=z_1^{n_1}z_2^{n_2}\ldots$, где
$z=(z_1,z_2,\ldots)\in{\Bbb{C}}^\infty$,
$n=p_1^{n_1}p_2^{n_2}\ldots$ -- каноническое разложение числа $n$
на простые множители, $\alpha(n)=(n_1,n_2,\ldots)$ -- финитный
мультииндекс. Следовательно, преобразование Лапласа для полугруппы
$\Bbb{N}^*$ имеет вид
$$
\widetilde k(\psi)=\sum_{n=1}^\infty
k(n)z^{\alpha(n)}.
$$
Достаточно доказать, что для любого полухарактера $\psi$
полугруппы $\Bbb{N}^*,\ |\psi|\leq \rho$
 найдется такое число   $s\in
\Bbb{C}, \ {\rm Re}s\geq 0$, что $|\widetilde k(\psi)-\widetilde
k(\rho\overline{\psi_s})|<C/2$, где $\psi_s(n):=n^{-s}$. Если
натуральное $m$ таково, что $\sum_{n=m+1}^\infty |k(n)|\rho(n)<C/8$, а
 $\psi(n)=\rho(n)z^{\alpha(n)}$, где
$z\in\overline{\Bbb{D}}^\infty$ (каждый полухарактер $\psi$
полугруппы $\Bbb{N}^*$, удовлетворяющий условию $|\psi|\leq \rho$,  имеет такой вид), то
$$
|\widetilde k(\psi)-\widetilde k(\rho\psi_s)|\leq\sum\limits_{n=1}^m
|k(n)|\rho(n)|z^{\alpha(n)}-n^{-s}|+C/4.
$$
 Выберем такие натуральные $r, M$, что любое натуральное
$n\leq m$ имеет каноническое разложение с простыми  множителями
$p_1,\ldots ,p_r$ и показателями $n_i\leq M$. Доказательство леммы
будет завершено, если мы убедимся, что найдется такое $s\in
\Bbb{C}, \ {\rm Re}s\geq 0$, для которого при всех $n_i\leq M$
справедливо неравенство
$$
|z_1^{n_1}\cdots z_r^{n_r}-p_1^{-n_1s}\cdots
p_r^{-n_rs}|<C/(4\|k\|_\rho).
$$
 Ввиду очевидного неравенства $|a^l-b^l|\leq l|a-b|$,
справедливого для $a,b\in\overline{\Bbb{D}}$, достаточно показать,
что за счет выбора $s\in \Bbb{C}, \ {\rm Re}s\geq 0$ можно сделать
сколь угодно малыми все разности $z_j-p_j^{-s}=z_j-e^{-s\log p_j}\
(j=1,\ldots ,r)$. Дело легко сводится к случаю $z_j=e^{-i\varphi_j},\
s=it$  с вещественными $\varphi_j,\ t$. Но поскольку числа $\log p_j$
линейно независимы над полем $\Bbb{Q}$, последнее утверждение
следует из аппроксимационной теоремы Кронекера (см., например,
\cite[(26.19, ii)]{HiR1}).

{\bf  Замечание}.  Аналог теоремы 3 в случае ненулевого
предела, вообще говоря, неверен, в чем нас убеждает следующий

\textbf{Пример 1 \cite{garman}.}  Пусть $S=\Bbb{N}^*$ --- мультипликативная полугруппа натуральных чисел, наделенная считающей мерой и дискретной топологией, $w=1$, и функция $q$ из
$L^1(\Bbb{N}^*)$ определена равенством $q(n)=1/2^{n+1}$. Так как
$|\sum_{n\geq 1}q(n)n^{-s}|\leq 1/2$ при ${\rm Re}s\geq 0$, то
для функции $k=u+q$, где  $u=1_{\{1\}}$ --- единица алгебры $L^1(\Bbb{N}^*)$, имеем
\[
\left|\sum\limits_{n\geq 1}\frac{k(n)}{n^s}\right|\geq 1-\left|\sum_{n\geq 1}\frac{q(n)}{n^{s}}\right|\geq \frac{1}{2} \quad(s\in
\Bbb{C},  {\rm Re}s\geq 0).
\]
Лемма 8 показывает теперь, что  $\widetilde q(\psi)\ne -1$  при всех
$\psi\in\; \widetilde S_1$. С помощью теоремы о неявной функции, примененной к банаховой алгебре
$L^1(\Bbb{N}^*)$,
выводим отсюда, что уравнение $K+q=-K\ast q$ имеет решение $K\in
L^1(\Bbb{N}^*)$. Так как $(u+K)(u+q)=u$, то $(1+\widetilde K(\psi)(1+\widetilde q(\psi))=1$
при всех $\psi\in\widetilde S_1$, а потому  $\widetilde K(\psi)\ne -1$. Если
мы положим $f\!=\!1\!+\!1\!\ast\! q$, то  $f\!\in\!
L^\infty(\Bbb{N}^*)$, причем $f\!+\!f\!\ast\! K\!=\!1$. С другой
стороны, $f(n)=1+\sum_{d|n}1/2^{d+1}\to 5/4$, если $n$ простое,
$n\to\infty$, и $f(p!)\geq 1+\sum_{1\leq j\leq p}1/2^{j+1}\geq
11/8$, для любого простого $p\geq 3$. Поэтому $f$ не имеет предела
на бесконечности.

В то же время, при дополнительных ограничениях на полугруппу   аналог теоремы 3 в случае ненулевого предела все же имеет место (см. \cite{Trudy}).

В заключение отметим следующее приложение полученных результатов.

{\bf Теорема 5.} {\it  Пусть $S$ -- свободная абелева
полугруппа с единицей $e$. Предположим, что на $S$ задана {\it
норма},\index{норма} т. е. такой гомоморфизм $N:S\to\Bbb{R}_+^*$,
что $N(a)>1$ при $a\ne e$. Обозначим через $P$ множество простых
элементов полугруппы $S$ (множество образующих). Тогда ряд
$\sum_{p\in P}1/N(p)$ расходится вместе с рядом }$\sum_{s\in
S}1/N(s)$.

Доказательство см. в \cite{garman}.

{\bf Пример 2.} Пусть $K$ --- поле алгебраических чисел
(конечное расширение поля $\Bbb{Q}$), $S$ --- полугруппа целых
идеалов этого поля, отличных от нуля. Это свободная полугруппа с
множеством образующих $\frak{P}$, состоящим из простых идеалов
$\frak{p_1}, \frak{p_2}, \ldots$, наделенная нормой $N$ (см.,
например, \cite{Leng}). Известно, что ряд $\sum_{\frak{a}\in
S}1/N\frak{a}$ расходится (точка $s=1$ является полюсом
дзета-функции Дедекинда $\zeta_K(s)=\sum_{\frak{a}\in
S}1/N\frak{a}^s$ поля $K$), поэтому расходится и ряд
$\sum_{\frak{p}\in \frak{P}}1/N\frak{p}$.

Замечание. Данная работа опубликована в \cite{taub}.

\end{document}